\documentclass[12pt] {amsart}

\newtheorem{theorem}{Theorem}[section]
\newtheorem{lemma}[theorem]{Lemma}
\newtheorem{proposition}[theorem]{Proposition}
\newtheorem{corollary}[theorem]{Corollary}

\theoremstyle{definition}
\newtheorem{definition}[theorem]{Definition}

\theoremstyle{remark}
\newtheorem{remark}[theorem]{Remark}

\numberwithin{equation}{section}

\def \real{{\mathbb R}}

\def \integer{{\mathbb Z}}

\def \Id {{\text{Id}\, }}
\def \Leb {{\text{Leb} }}

\def\BB{{\mathcal B}}
\def\CC{{\mathcal C}}
\def\DD{{\mathcal D}}

\def\FF{{\mathcal F}}

\def\II{{\mathcal I}}

\def\LL{{\mathcal L}}
\def\MM{{\mathcal M}}

\def\SS{{\mathcal S}}

\def\XX{{\mathcal X}}

\def\today {\ifcase\month\or January \or February \or March \or
April \or May \or June
\or July \or August \or September \or October \or November \or December
\fi
\number\day~\number\year}

\begin{document}
\title[Anisotropic Sobolev spaces]
{Anisotropic Sobolev spaces and dynamical transfer
operators: $C^\infty$ foliations}

\author{Viviane Baladi}
\address{CNRS-UMR 7586, 
Institut de Math\'e\-ma\-ti\-ques  Jussieu, Paris, France}

\date{December 2004}

\begin{abstract}
We consider a $C^\infty$ Anosov diffeomorphism $T$ with
a $C^\infty$ stable dynamical foliation. We show upper bounds
on the essential
spectral radius of its transfer
operator acting on  anisotropic Sobolev spaces. 
(Such bounds are related   to the
essential decorrelation rate for the SRB measure.)
We compare our results to the estimates of Kitaev
on the domain of holomorphy of dynamical determinants
for differentiable dynamics.
\end{abstract}


\thanks{Thanks to A. Avila, L. Boutet de Monvel,
G. David, P. G\'erard, S. Gou\"ezel, F.~ Ledrappier,
C. Liverani, and M. Tsujii for useful conversations. }

\maketitle


\section{Introduction}

Let $T$ be an Anosov diffeomorphism on a 
$d$-dimensional compact connected $C^\infty$ Riemann manifold
$\XX$ (i.e., $T\XX =E^u \oplus E^s$ and there are
$C>0$, $\gamma >1$,
with $|DT ^n|E^s|\le C \gamma^{-n}$,
$|DT ^{-n}|E^u|\le C \gamma^{-n}$ for all
$n\ge 1$). Denote the Jacobian of $T$ with respect to
Lebesgue by $|\det DT|$.
To construct SRB measures and to analyse their
speed of mixing, it is natural to consider 
the following operators, defined initially
on $C^\infty$ functions:
\begin{equation}
\label{xfer}
\MM \varphi=\frac{ \varphi\circ T^{-1}}{|\det DT|\circ T^{-1}}\, ,
\qquad
\LL \varphi = \varphi \circ T \, .
\end{equation}
The operator $\LL$ fixes the constant functions, while $\MM$ fixes
the constant functions if and only if $\det DT$ is constant (i.e., if
$T$ is volume preserving).  The dual of $\MM$ restricted
to elements of the dual of $\CC^\infty$ which are finite complex
measures, absolutely continuous
with respect to Lebesgue with a $\CC^\infty$ density, coincides  with
$\LL$ viewed as acting on the density and vice-versa. Alternatively,
the dual of $\MM$ acting on $L^1(\XX, \Leb)$ is $\LL$ acting on 
$L^\infty(\XX, \Leb)$.

For $w \in \XX$, and $\widetilde T$ an
Anosov diffeomorphism on $\XX$, introduce local hyperbolicity exponents
($|\cdot|$ denotes euclidean norm)
\begin{equation}\nonumber
\begin{split}
\lambda_{w}(\widetilde T)^{-1}
&=\sup_{v\in E^u(\widetilde T(w)), |v|=1} |D_{T(w)} \widetilde T^{-1} (v)|\, ,\cr
\nu_{w}(\widetilde T)&=\sup_{v\in E^s(w), |v|=1} |D_w\widetilde T (v)|\, .
\end{split}
\end{equation}
Assume $T$ is $C^{r+1}$ for some $r>0$.
Kitaev \cite{Ki} proved that the following ``dynamical Fredholm determinant''
\begin{equation}\nonumber
\nonumber
d(z)=\exp -\sum_{n=1}^\infty
\frac {z^n} {n} \sum_{T^n(x)=x} \frac {1}{ |\det (DT^n(x)-\Id)]}
\end{equation}
extends to a holomorphic function in each disc
$\{z\mid |z| \cdot \rho^{(p,s)}_1(T)<1\}$, where  $p\in (-r,0)$,
$s=r+p$, and
\begin{equation}
\nonumber
\rho^{(p,s)}_1(T)=\lim_{n\to \infty}
\biggl ( \int_\XX
\max\bigl ( (\lambda_{w} (T^n))^{p}, (\nu_{w} (T^n))^{s}
\bigr )  \, dLeb(w) \biggr )^{1/n} <1 \, .
\end{equation}
One may take $s=-p=r/2$: Kitaev's result is then reminiscent of
the ``loss of one half of the H\"older exponent'' which occurs when going from
two-sided subshifts to one-sided subshifts in symbolic dynamics \cite{Bow},
since one easily sees that $\rho^{(-r/2,r/2)}_1(T)\le \gamma^{-r/2}$.

In view of the results of Ruelle \cite{Ru} for smooth expanding maps,
it is natural to look for  Banach spaces $\BB_{p,s,\LL}$, respectively
$\BB_{p,s,\MM}$,  on which the essential spectral
radius of $\LL$, respectively $\MM$, is $\le \rho^{(p,s)}_1$.
Set
\begin{equation}\nonumber
\rho^{(p,s)}_\infty(T)=\lim_{n\to \infty}
\biggl ( \sup_\XX
\max\bigl ( (\lambda_{w} (T^n))^{p}, (\nu_{w} (T^n))^{s}
\bigr )   \biggr )^{1/n} <1 \, .
\end{equation}
Clearly $\rho^{(p,s)}_\infty(T)\ge \rho^{(p,s)}_1(T)$ and, e.g.,
 $\rho^{(-r/2,r/2)}_\infty(T)\le \gamma^{-r/2}$.

\medskip

We shall assume that  $T$ is $C^\infty$ and the stable foliation of
$T$ (or its unstable foliation) is $C^\infty$.  
(This is a very strong assumption, and the corresponding
case should be viewed as a ``toy model" in which the features of
our symbolic calculus  approach are completely transparent:
The heart of the proof is contained in a half page,
between (\ref{heart2}) and (\ref{heartn}) below.) 
We introduce in Subsection ~\ref{space},
for $p$,  $s$  in $\real$ and   $1<t< \infty$, a Banach
space $W^{p,s-p,t}(\XX)=W^{p,s-p,t}(\XX, T)$ of distributions, 
based on $L^t(\Leb)$.
\footnote{Controlling the spectrum on a scale of Sobolev spaces may 
be useful: see \cite{BB}.}

\smallskip

Our main result (Theorem \ref{mainthm}) when the stable foliation is $C^\infty$
is that, if $T$ is volume preserving,
the essential spectral radius $\rho_{ess}$ of $\LL$ on $W^{p,s-p,t}(\XX)$ is at most $\rho^{(p,s)}_\infty(T)$
for all $p <0$, $s >0$ and
$t \in (1, \infty)$; while if $T$ does not preserve
volume thn $\limsup_{t\to \infty}\rho_{ess}(\LL|_{W^{p,s-p,t}(\XX)})
\le \rho^{(p,s)}_\infty(T)$. 
If the unstable foliation is $C^\infty$,  
the essential spectral radius
of $\MM$  on $W^{p,s-p,t}(\XX, T^{-1})$ is at most 
$\lim_{n \to \infty} \sup_\XX |\det DT^n|^{-(t-1)/tn}  \cdot \rho^{(-s,-p)}_\infty(T)$
for all $p <0$,  $s>0$ and $t\in (1, \infty)$ (Theorem~\ref{mainthm2}).

\smallskip

Propositions ~\ref{L1} and
~\ref{L1,2} give upper bounds related to
$\rho^{(p,s)}_1(T)$: They imply 
\begin{equation}\nonumber
\begin{split} 
\limsup_{t\to \infty}&\rho_{ess} (\LL|_{W^{p,s-p,t}})\le
\lim_{n\to \infty} \| \det DT^n|_{E^u} \|_{L^\infty(\Leb)}^{1/n} \rho^{(p,s)}_1(T)\, ,\cr  
\limsup_{t\to 1}
&\rho_{ess} (\LL|_{W^{p,s-p,t}})\le \lim_{n\to \infty} \| (\det DT^n|_{E^s})^{-1} \|_{L^\infty(\Leb)}^{1/n} \rho^{(p,s)}_1(T)
\, .
\end{split}
\end{equation}

\medskip

Finally, we study in the appendix  the essential spectral radii of
$$\LL_t \varphi =|\det DT|^{1/t} \cdot (\varphi \circ T) \, ,
\qquad \MM_t\varphi =\frac {\varphi \circ T^{-1} } {|\det DT|^{1-1/t} \circ T^{-1}}\, .$$

\medskip

The case when $T$ is $C^{r+1}$ (for some $r>0$)  and neither of the
dynamical foliations is $C^\infty$, but at least one of them
is $C^{1+\epsilon}$ (for $\epsilon>0$) will be treated in a forthcoming work
\cite{Ba},
using spaces due to Alinhac \cite{A2}. We hope that the (general)
$C^\alpha$ foliation case will be amenable to the present
approach.  Gou\"ezel and Liverani \cite{GL} have independently obtained
non trivial, but weaker, bounds for the essential spectral
radius of $\MM$, on a different Banach space,  
in this general case.

\smallskip

We end this introduction with three open problems:

\begin{remark}[Links with SRB measures]
With the techniques of Blank--Keller--Liverani \cite{BKL}, one should  
obtain that the spectral radius of $\LL$
on each $W^{p,q,t}(\XX)$ is one, that $1$ is
is a semi-simple eigenvalue, and that the corresponding
eigenvector (in the dual   of $W^{p,q,t}(\XX)$)
for the dual of $\LL$ is  an invariant probability measure $\mu$
with ergodic basin of full Lebesgue measure.
Furthermore, the multiplicity of the eigenvalue $1$ is equal to the
number of ergodic components of $\mu$, and each
ergodic component is an SRB measure. Also, if $1$
is a simple eigenvalue then it  is the only eigenvalue
on the unit circle: this corresponds to exponential decay of correlations
for smooth observables. 
If the unstable foliation is $C^\infty$,
the SRB measure(s) of $T$ can alternatively be constructed
with the fixed point of $\MM$ in $W^{p,q,t}(\XX, T^{-1})$.
\footnote{The operators $\LL_t$, $\MM_t$  ``interpolate'' between
the SRB measures of $T$, $T^{-1}$.}
\end{remark}

\begin{remark}[Spectral stability]
The perturbation techniques of \cite{BKL},  should imply
stability  of the spectrum of $\LL$ (including eigenprojectors) outside
a disc of radius $\rho$,  under stochastic
and $C^{r+1}$ deterministic perturbations of $T$,
perhaps up to  taking  $\rho >\rho^{(p,s)}_\infty$. 
For deterministic perturbations $\widetilde T$ of $T$,
the Banach spaces $W^{p,q,t}(\XX,T)$
and $\widetilde W^{p,q,t}=W^{p,q,t}(\XX,\widetilde T)$ 
are
different. ``Stability of the
eigenprojector'' $\Pi$ of $T$ associated to an eigenvalue $\tau$ of large
enough modulus means the following 
(assume $\tau$ is simple): Let $\widetilde \LL$ denote the transfer operator of $\widetilde T$; then, 
if $\widetilde T$ is close enough to
$T$, there are a Banach space $ W_\epsilon$ contained in the intersection
of $\widetilde W^{p,q,t}$ and $W^{p,q,t}$, a rank-one projector  
$ \Pi_\epsilon$
(on $W_\epsilon$, $\widetilde W^{p,q,t}$, and $W^{p,q,t}$), and
 a simple eigenvalue $(\tilde \tau,\widetilde \Pi)$ for $\widetilde \LL$
on $\widetilde W^{p,q,t}$,
so that both $\|\Pi_\epsilon-\Pi\|_W$ and 
$\| \Pi_\epsilon-\widetilde \Pi\|_{\widetilde W}$ are small.
\end{remark}

\begin{remark}[Essential decorrelation radius]
For $C^{r+1}$ expanding circle endomorphisms 
$F$, the essential spectral
radius $\rho_{ess}(\MM_F|_{C^r})$
of $\MM_F \varphi(x)=\sum_{F(y)=x} \varphi(y)/|\det DF(y)|$ acting
on $C^r$ functions (see \cite{CE} and references therein)
is equal to
$$\lim_{n\to \infty} \bigl (\int |(F^n)'(x)|^{-r} \,  dLeb(x) \bigr )^{1/n}=
\lim_{n\to \infty} \bigl (\int |(F^n)'(x)|^{-r} \,  d\mu_{SRB}(x)\bigr  )^{1/n} \, . $$ 
However,
for $C^{r+1}$ expanding maps in arbitrary dimension 
\cite{GLa}
\begin{equation}\label{GuL}
\begin{split}
\rho_{ess}(\MM_F|_{C^r})=&\exp (\sup_{\mu}  
\{ h_\mu - \int \log |\det DF | \, d\mu - r \cdot \chi_\mu\}) \cr
&\qquad \qquad\le
\lim_{n\to \infty} \bigl (\int \sup_{|v|=1} |D_x(F^n)(v)|^{-r} \,  dLeb(x)\bigr )^{1/n} \, ,
\end{split}
\end{equation}
where $\mu$ ranges over ergodic $F$-invariant probability measures,
$h_\mu$ is the entropy of $\mu$,
and $\chi_\mu$ denotes the smallest (positive) Lyapunov
exponent of $DF$.
The inequality in (\ref{GuL}) can be strict.
In the other direction, note that
 $\rho_{ess}(\MM_F|_{C^r})\ge  \exp(-r \chi_{\mu_{SRB}})$,
and  the inequality can be strict \cite{CE}, even in dimension one.
The results of Avila et al.  \cite{AGT},
indicate that in dimension at least two there may be
Banach spaces containing all $C^r$ functions
on which the essential spectral radius of $\MM_F$
is strictly smaller than $\rho_{ess}(\MM_F|_{C^r})$.
(This would imply  \cite{CE} that $\rho_{ess}(\MM_F|_{C^r})$ 
may be strictly larger than the essential decorrelation
radius of $F$ for $C^r$ observables and thus
$\rho_{point-ess}(\MM_F|_{C^r})<\rho_{ess}(\MM_F|_{C^r})$.)

Let $T$ be a transitive $C^{\infty}$ Anosov diffeomorphism with
both foliations $C^\infty$.
Let $\rho_{ess}^+(p,s,t)$ and  $\rho_{ess}^-(p,s,t)$
be  the essential spectral radii  of $\LL$ acting
on $W^{p,s-p,t}(\XX)$ and $W^{-p,-s+p,t}(\XX, T^{-1})$, respectively, and set
$$ \rho(r):=
\min\bigl (\inf_{\stackrel{t, p \in(-r,0)}{  s\in (0,r+p)}} \rho^+_{ess} (p,s,t),
\inf_{\stackrel{t, p \in(-r,0)} { s\in (0,r+p) }} \rho^-_{ess} (p,s,t) \bigr )
\, . 
$$
We expect that
$\inf_{\BB_r} \rho_{ess}(\LL|_{\BB_r})$, where $\BB_r$
spans all  Banach spaces  of distributions of order $\le r$,
containing all $C^r$ functions, and  on which $\LL$ acts boundedly,
coincides with
the essential decorrelation radius $\hat \rho(r)$ of $T$ for $C^r$
functions, and that
$\hat \rho(r) < \rho(r)$
can occur. 
\end{remark}

\section{Bounding the essential spectral radius}

\subsection{Preliminaries} \label{prel}
From now on and until the end of Subsection~ \ref{proofLY}, $T$ is Anosov and
$C^\infty$, with a $C^\infty$ stable foliation $\FF^s$.
Write $I=(-1,1)$, and let $d_s$ be the dimension of $\FF^s$.
We work with  $C^\infty$ foliated charts $\kappa$, $V$:
let $\cup_{i \in \II} V_i$ be a finite covering of
$\XX$ by small open sets, and let 
$U_i=I^d=I^{d_s}\times I^{d-d_s}$  be $\#\II$ disjoint copies
of $I^d$, viewed as subsets of disjoint copies
$\real^d_i$ of $\real^d$. 
Let $\kappa_i : V_i \to U_i$
be $C^\infty$ diffeomorphisms so that
$\kappa_i^{-1}$ of each horizontal segment
is the intersection of a leaf of $\FF^s$ with $V_i$.
In addition, we require
that $\kappa_i^{-1}(\{(0,y)\})$ is the unstable leaf of $\kappa_i^{-1}(0,0)$
intersected with $V_i$ (this is a way to require  closeness 
of the vertical foliation in $I^d$ and the image of leaves of the unstable foliation
$\FF^u$).

Choose a $C^\infty$ partition of the unity
$\{\psi_i\}$ on $\XX$, compatible with
$V=\{V_i\}$, i.e., each $\psi_i$ is supported in $V_i$.
Then, for each $n\ge 1$
\begin{equation}
\label{maintransfer}
\LL ^n\varphi (w)=\sum_{i, j} \psi_j(T^n (w)){\psi}_{i}(w)\cdot \varphi (T^n(w))  \, .
\end{equation}

If $V_{ij}=V_{ij,n}:=T^{-n}(V_j)\cap V_i \ne \emptyset$, setting
$U_{ij,n}:=\kappa_i(T^{-n}(V_j)\cap V_i)\subset U_i$,
the map  $T^n_{ij}:U_{ij,n}\to U_j$
has a derivative in block form:
\begin{equation}\nonumber
\left ( \begin{array}{cc}
A^{tr}_{ij,n}(x,y) & B^{tr}_{ij,n}(x,y) \\
0 & D^{tr}_{ij,n}(x,y) \\
\end{array} \right ) \, , \, (x,y)\in (I^{d_s}, I^{d-d_s}) \, ,
\end{equation}
with  $A_{ij,n}$ a $d_s\times d_s$ matrix, $D_{ij,n}$ a $(d-d_s)\times (d-d_s)$ matrix, and
\begin{equation}\label{rem0}
\begin{split}
&| A_{ij,n}(x,y)|\le\nu_{ij}(T^n) := \sup_{w\in V_{ij}} \nu_w(T^n)<1\, ,\cr
& | D_{ij,n}(x,y) ^{-1}|\le\lambda_{ij}(T^n)^{-1}:= 
\sup_{w\in V_{ij}} (\lambda_{T^n(w)}(T^n) )^{-1}<1\ . \cr
\end{split}
\end{equation}
Furthermore, for each $\epsilon$, there exists $\delta$ so that if 
${\rm diam}\, V <\delta$ then 
\begin{equation}
\label{rem}
|B_{ij, n}(x,y) v |\le \epsilon  | D_{ij, n}(x,y) v | \, , \forall n\ge 1\, , \, 
\forall v \in \real^{d-d_s} \, .
\end{equation}

\subsection{Elementary spaces $W^{p,q,t}(\real^d)$}
Let $p$ and $q$ be real numbers.
We introduce the ``symbol'' $a_{p,q}(\xi,\eta)$, for
$(\xi, \eta) \in \real^{d_s}\times \real^{d-d_s}$:
$$
a_{p,q}(\xi,\eta)= 
(1+|\xi|^2 +|\eta|^2)^{p/2} (1+|\xi|^2 )^ {q /2}\, . 
$$
The corresponding linear operator $a^{Op}_{p,q}$ maps the
space $\SS$ of rapidly decaying $C^\infty$ functions on $\real^d$
into itself via
$$
a^{Op}_{p,q} (f) (x,y)
=(2\pi)^ {-d}\int \int e^{ix\xi} e^{iy\eta} a_{p,q}(\xi,\eta) \hat f(\xi,\eta)\,  d\xi \, d\eta \, ,
$$
where  the Fourier transform of $f$ is
$
\hat f(\xi,\eta)=\int \int e^{-ix\xi} e^{-iy\eta} f(x,y) \, dx \, dy 
$.

\begin{definition}[Anisotropic Sobolev spaces]
For $1\le t\le\infty$, let
$W^{p,q,t}(\real^d)$   be the closure
of $\{ f \in \SS(\real^d)\}$ for the  $L^t(\real^d, \Leb)$ norm of
$a^{Op}_{p,q} (f)$, with induced norm,
denoted $\|\cdot\|_{p,q,t,\real^d}$.  
\end{definition}

By construction, $a^{Op}_{p,q}$ extends to a
bounded invertible operator  from $W^{p,q,t}(\real^d)$ to $L^t(\real^d)$.
Clearly, $H^{p,q}(\real^d)=W^{p,q,2}(\real^d)$  is  a Hilbert space.

\begin{lemma}[Boundedness/compactness of embedding]\label{compact}
Assume that $1<t<\infty$. Denote by
$W^{p',q',t}_C(\real^d)$ those $f\in W^{p',q',t}(\real^d)$ 
supported in a compact subset of $\real^d$.
If $q' \ge q$ and $p'\ge p$ then the natural injection
$W^{p',q',t}_C(\real^d) \subset W^{p,q,t}(\real^d)$
is bounded. This
injection is  compact if $q' \ge q$ and $p'>p $ .
\end{lemma}

\begin{proof}
If $t=2$, the proofs of Theorems 2.5.2 and 2.5.3 in \cite{H0} adapt readily.
The general case is an easy exercise.
\end{proof}

\begin{remark}
More generally, we may introduce
classes  of (symbols) of pseudodifferential
operators:
Let $p$ and $q$ be real numbers.
We say that  $b\in C^{\infty}(I^d \times \real^d, \real)$, 
belongs to $S^{p, q}$ if for any
multi-indices\footnote{We decompose multi-indices $\gamma=(\gamma', \gamma'')$
in this way tacitly from now on.}  $\alpha=(\alpha',\alpha'')$ and
$\beta=(\beta',\beta'')$
in $\integer^{d_s+(d-d_s)}_+$,  there exists $C_{\alpha, \beta}$
so that
\begin{equation*}
\begin{split}
\sup \left |  
\partial^{\alpha'}_\xi \partial^{\alpha''}_\eta \partial^{\beta'}_x
\partial^{\beta''} _y b(x,y,\xi,\eta)\right |
&\le  C_{\alpha, \beta} (1+|\xi|+|\eta|)^{p-|\alpha''|} (1+|\xi|)^{q-|\alpha'|} \, .
\end{split}
\end{equation*}
The spaces $S^{p,q}$ and $H^{p,q}$ were studied by
Kordyukov \cite{Ko}.
The 1963  edition of H\"ormander's book \cite[II.2.5]{H0} contains a treatment
of a special case of the spaces $S^{p,q}$ when $d_s=1$.
See  also Sabl\'e-Tougeron \cite{ST} for applications of these 
special cases.
\end{remark}

\subsection{Banach spaces $W^{p,q,t}(\XX)$
and Leibniz formula}\label{space}
Let $\kappa$, $V$ be a chart 
and  $\psi$ be a compatible partition
of unity as in Subsection~ \ref{prel}.

\begin{definition}\label{defnorm}
Let $p$, $q$ be real numbers, and let $1\le t\le \infty$.
$W^{p,q,t}(\XX, \kappa,V,\psi)$ is 
$\{ \varphi \in \DD'(\XX)
\mid (\psi_i  \cdot \varphi)\circ \kappa_i^{-1} \in W^{p,q,t}(\real^d_i)\, ,
\forall i \in \II \}$,
normed by
$$
\|\varphi\|_{p,q,t}=\sum_{i\in \II} \|(\psi_i \cdot \varphi)\circ \kappa_i^{-1}\|_{p,q,t,\real^d_i} \, .
$$
\end{definition}

\begin{remark} 
If $1<t< \infty$, the Banach spaces $W^{p,q,t}(\XX,\kappa,V,\psi)$
are independent of the charts $(\kappa,V)$ and of the
partition of unity $\psi$: A version of the change of variables
theorem for pseudodifferential operators, see e.g. \cite[I.7.1]{AG},  shows
that the norms corresponding to different $(\kappa, V, \psi)$
are equivalent. (See
Lemmas ~\ref{Leibniz} and~ \ref{comp} below.)
We may thus write $W^{p,q,t}(\XX)$.
$H^{p,q}(\XX)=W^{p,q,2}(\XX)$ is a Hilbert space.
\end{remark}

\begin{remark}
$W^{p,q,t}(\XX)$  is the Banach space of
distributions $f$  on $\XX$ so that $(1+\Delta_s)^{q/2} (1+\Delta)^{p/2} f
 \in L^t(\XX)$, with the induced $L^t(\XX)$ norm,
where $\Delta$ is the Laplacian and $\Delta_s$ is the stable foliated  Laplacian.
In particular, if $p\le 0$ and $0\le q  \le  r$, it contains all $C^r$ functions.
\end{remark}

We start with a useful remark:

\begin{lemma}[Proper support]\label{proper}
For every compact subsets $K_1$, $K_2$ of $I^d$, with
$K_1$ included in the interior of $K_2$, and
for every $1<t<\infty$,  $p$, $q$,  there are $C>0$ and
a $C^\infty$ function $\Psi_K:\real^d\to [0,1]$, supported in
$K_2$, so that
for
each $f \in W^{p,q,t}(\real^d)$ supported in $K_1$,
$$
\|\Psi_K \cdot a_{p,q}^{Op}(f)-a_{p,q}^{Op}(f)\|_{L^t} \le C
\| f\|_{p-1,q,t} \, .
$$
\end{lemma}

\begin{proof}
Using that the kernel of a pseudo-differential operator
is $C^\infty$ outside of the diagonal, a standard construction allows
to write $\Psi\cdot a-a$ (acting on compactly supported distributions)
as an operator with a $C^\infty$ kernel
(see e.g. \cite[Prop 6.3]{AG}). Integrate by parts to conclude.
\end{proof}

\begin{lemma}[Leibniz formula]\label{Leibniz}
Let $1<t<\infty$, let $p$, $q$ be real numbers and
let  $h$ be a compactly supported
$C^\infty$ function on $I^d$. 
Then there is $C(h)>0$,
and there exists $C_x(h)$, depending only on
\begin{equation}
\label{opM}
\sup_{|\beta'|\in \{1, 2\}, (x,y)\in I^d} |\partial^{\beta'}_x h(x,y)|\, , 
\end{equation} so that for every
$f \in W^{p,q,t}(\real^d)$
\begin{equation}\label{Leibnizz}
a^{Op}_{p,q}( h \cdot f) =h \cdot a^{Op}_{p,q}(  f)+
g_1 + g_2  \, ,
\end{equation}
with
$\|g_1\|_{L^t} \le C_x(h)  \| f\|_{p,q-1,t}$
and $\|g_2\|_{L^t} \le C(h)  \| f\|_{p-1,q,t}$.
\end{lemma}

\begin{proof} Multiplication by $h$ is a pseudodifferential operator.
Composing it with $a^{Op}_{p,q}$, we get a new operator 
$b^{Op}$.
Using
a  Taylor series of order one
(see e.g. \cite[Th\'eor\`eme I.4.1 and \S I.8.2]{AG}), we find $b (x,y,\xi,\eta)=$
\begin{equation}\nonumber
\begin{split}
& a_{p,q}(\xi,\eta) \cdot h(x,y)+  \frac {2}{ (2\pi)^{d} }\sum_{|\alpha|+|\beta|=2} 
\frac {(-1)^{|\alpha|+|\beta|} }{ \alpha ! \beta !} \int_0^1 (1-s) \cdot \cr
&\qquad\, \cdot  \int e^{-i (u,v)(\omega,\theta)} 
\omega^{\beta'} \theta^{\beta''} u^{\alpha'} v^{\alpha''}
 \cdot \partial^{\beta'}_\omega \partial^{\beta''}_\theta
a_{p,q}(\xi-s\omega, \eta-s\theta)  \, d\omega\, d\theta\cr
&\qquad\qquad\qquad \qquad \qquad\qquad \cdot 
\partial^{\alpha'}_u\partial^{\alpha''}_v h(x-s u , y-s v) 
 \, du \, dv \, ds  \, .\cr
\end{split}
\end{equation}
The 
symbol  $a_{p,q}(\xi,\eta) \cdot h(x,y)$ gives rise to
the first term in the right-hand-side of (\ref{Leibnizz}). 
For the remainder term, the usual integrations by parts
\cite[\S I.8.2, p.56]{AG} yield
a linear combination of terms
\begin{equation}
\nonumber
\begin{split}
b ^{\gamma, j}(x,y,\xi,\eta)&:= \int_0^1 (1-s)  s^j \cdot  \int e^{-i (u,v)(\omega,\theta)} \cdot 
\partial^{\gamma'}_\omega\partial^{\gamma''}_\theta 
a_{p,q}(\xi-s\omega, \eta-s\theta)  \cr
&\qquad \qquad\quad \cdot
 \partial^{\gamma'}_u \partial^{\gamma''}_v h(x-s u , y-s v) 
\, d\omega d\theta \, du dv \, ds  \, ,\cr
\end{split}
\end{equation}
where $j\in \{0,1,2\}$ and $|\gamma'|+|\gamma''|\in\{ 1,2\}$
(the number of terms and the coefficients in the linear
combination are independent of $h$ and $a_{p,q}$). 
If $|\gamma''|=0$ then $|\gamma'|\in\{1,2\}$,
and this gives   $g_1$, as we explain next.
Define a symbol $\tilde b=b ^{\gamma, j} (a_{p,q-1})^{-1}$.
By \cite[Th\'eor\`eme~9]{CM} it 
suffices to show that there is
$C_x(h)$ so that $\sup | \partial^{ \alpha}_{x,y} 
\partial^{\beta}_{\xi,\eta} \tilde b(x,y,\xi,\eta)|\le C_x(h)(1+|\xi|+|\eta|)^{-|\beta|}$, 
for all $| \alpha|\le 1$ and all $ \beta$.
This can be shown by a straightforward (although tedious) implementation
of the standard oscillatory integral argument \cite[\S I.8.2, p.56]{AG}. 
Finally, if $|\gamma''|\ge 1$ then the  term corresponding
to  $b ^{\gamma, j}$ may be included
in  $g_2$, working with $b ^{\gamma, j} ( a_{p-1,q})^{-1}$.
\end{proof}

\subsection{Bounding the essential spectral radius of $\LL$}

\begin{theorem}\label{mainthm}
Let $T$ be a $C^\infty$ Anosov diffeomorphism on a compact 
manifold, with a $C^\infty$ stable foliation.   For any $p<0$, $s>0$,
and $t\in(1,\infty)$,
the essential spectral radius of $\LL$  on $W^{p,s-p,t}(\XX)$  is not larger than
$\lim_{n \to \infty} \sup_\XX |\det DT^n|^{-1/tn}  \cdot \rho^{(p,s)}_\infty(T)$.
\end{theorem}

Note that the essential
spectral radius of the dual of $\LL$
i.e. (an extension of) 
$\MM$ acting on the dual  of $W^{p,s-p,t}(\XX)$ coincides with
the essential spectral radius of $\LL$ on $W^{p,s-p,t}(\XX)$.  Also, if the unstable foliation 
is $C^\infty$, then   $\varphi \mapsto \varphi \circ T^{-1}$ 
on  $W^{p,s-p,t}(\XX,T^{-1})$ 
has essential spectral radius $\le
\lim_{n \to \infty} \sup_\XX |\det DT^{-n}|^{-1/tn}  \cdot\rho^{(p,s)}_\infty(T^{-1})$
(note that $\rho^{(p,s)}_\infty(T^{-1})=\rho^{(-s,-p)}_\infty(T)$).

It is convenient to extend each
$T^n_{ij}=
\kappa_j\circ T^n \circ\kappa_i^{-1}:U_{ij,n} \to U_j$ to a $C^\infty$ diffeomorphism
from $\real^d_i$ onto its image in $\real^d_j$ in such a way
that the intersection of $U_j$
with the image of $U_i$ by the extended map coincides with $T^n_{ij}(U_j)$,
and so that the extended map is the identity outside of a large compact set.
(The extension is still noted $T^n_{ij}$.)
The theorem will  be a consequence of the following lemma,
proved in \S \ref{proofLY}: 

\begin{lemma}[Lasota-Yorke inequality]\label{comp}
There exist $\delta_0 >0$ and $C_0$, so that for each cover  with
${\rm diam} \, V <\delta_0$, 
and for each  $n\ge 1$  there exists $C(n)>1$,
so that for  every 
$f \in W^{p,q,t}(\real^d_j)$, compactly supported
in $U_j$, and  each $C^\infty$ function $\Psi_{ij}:\real_i^d
\to [0,1]$ compactly
supported in $U_{ij,n}$ 
\begin{equation}\nonumber
\begin{split}
\|   & \Psi_{ij} \cdot a_{p,q}^{Op}(f \circ   T_{ij}^n)\|_{L^t}
 \cdot \inf_{V_{ij}} |\det DT^n|^{1/t}\cr
&\qquad \le C_0  \cdot
\max((\lambda_{ij}(T^n))^{p},(\nu_{ij}(T^n))^{q+p}) 
\|  f \|_{p,q,t,\real^d_j}\cr
&\qquad\qquad\quad+ C(n) \|  f \|_{p-1/2,q,t,\real^d_j} \, , \, \forall
p\le 0\, , \, q\ge -p\, , \, 1<t<\infty \, .
\end{split}
\end{equation}
\end{lemma}

\begin{proof}[Proof of Theorem \ref{mainthm} using Lemma~\ref{comp}]
Let $\delta_0$ be as in Lemma~\ref{comp}.
For $\delta <\delta_0$,  let
$(\kappa, V)$ be a foliated chart of diameter at most 
$\delta$, and let $\psi$ be an adapted partition of unity.
Set  $f_j|_{U_j}=(\psi_j \cdot \varphi) \circ \kappa_j^{-1}$, extending
by zero on $\real^d_j$. By definition, for all $n\ge 1$,
$$
\| \LL^n \varphi \|_{p,q,t}\le
\sum_i \sum_{j :T^n(V_i)\cap V_j \ne \emptyset}
\| (\psi_i \circ \kappa_i^{-1})  
\cdot (  f_j \circ T^n_{ij}) \|_{p,q,t,\real^d_i}\, .
$$
By Lemma \ref{proper}, 
there is a $C^\infty$ function
$\Psi_{ij}:\real^d_i\to [0,1]$, supported
in a compact subset of $U_{ij,n}$, so that $\| (\psi_i \circ \kappa_i^{-1})  
\cdot ( f_j \circ T^n_{ij})\|_{p,q,t,\real^d_i}\le$
\begin{equation}\nonumber
\begin{split}
&  \| \Psi_{ij} a_{p,q}^{Op}( (\psi_i \circ \kappa_i^{-1})
f_j \circ T^n_{ij} )\|_{L^t}
+C  \|(\psi_i \circ \kappa_i^{-1})
(f_j \circ T^n_{ij}) \|_{p-1,q,t} \, .
\end{split}
\end{equation}
By  Lemma \ref{Leibniz},  
the first term in the above sum is bounded by
\begin{equation}\nonumber
\begin{split}
& C(\psi) \cdot \| \Psi_{ij} a_{p,q}^{Op}( f_j \circ T^n_{ij} )\|_{L^t}
+  C(\psi) \cdot \| ( f_j \circ T^n_{ij})\|_{p-1,q,t} \, .
\end{split}
\end{equation}
Set $\rho(p,s,n)=\max_{i,j} \max(( \lambda_{ij}(T^n))^p,(\nu_{ij}(T^n))^s))$.
By Lemma~\ref{comp} 
\begin{equation}\nonumber
\begin{split}
&\| \LL^n \varphi \|_{p,s-p,t}\cdot\inf_{V_{ij}} |\det DT^n|^{1/t}
\cr
&\qquad\qquad\le C_0 C(\psi) \rho(p,s,n) 
\cdot\sum_{i,j}  
\| a_{p,s-p}^{Op}( (\psi_j \cdot \varphi)\circ \kappa_j^{-1})\|_{L^t}\cr
&\qquad\qquad\qquad+  C(n, \psi)  \sum_{i,j}
\|(\psi_j \cdot \varphi)\circ \kappa_j^{-1}\|_{p-1/2,s-p,t}\cr
&\qquad\qquad\le\# \II \cdot  C_0 C(\psi)  \rho(p,s,n)       
 \cdot 
\sum_j \|(\psi_j \cdot \varphi) \circ \kappa_j^{-1}\|_{p,s-p,t}            \cr
&\qquad\qquad
\qquad+ \# \II \cdot C(n,\psi) 
\sum_j \| (\psi_j \cdot \varphi)\circ \kappa_j^{-1}\|_{p-1/2,s-p,t}\cr
&\qquad \qquad\le  C_1  
\rho(p,s,n) \cdot 
 \|  \varphi\|_{p,s-p,t}  + C_2(n)\|  \varphi\|_{p-1/2,s-p,t}  
\,.\cr
\end{split}
\end{equation}
By Lemma \ref{compact},  we can apply Hennion's theorem \cite{He}.
\end{proof}

\subsection{Bounds involving averaged hyperbolicity exponents}

\begin{proposition}\label{L1}
Let $T$ be a $C^\infty$ Anosov diffeomorphism on a compact 
manifold, with a $C^\infty$ stable foliation. 
For any $p<0$, $s>0$, and   $1<t< \infty$,
the essential spectral radius of $\LL$  on
$W^{p,s-p,t}(\XX)$  is  
\begin{equation}
\begin{split}\nonumber
&\le 
\lim_{n\to \infty}
\biggl ( \int_\XX
\max\bigl ( (\lambda_{w} (T^n))^{p}, (\nu_{w} (T^n))^{s}
\bigr ) \cr
&\qquad\qquad\qquad\qquad  \cdot | \det DT^n|_{E^u} | \cdot
|\det DT^n |^{-1/t} \, dLeb(w) \biggr )^{1/n} \cr
&
=
\lim_{n\to \infty}
\biggl ( \int_\XX
\max\bigl ( (\lambda_{w} (T^{n}))^{p}, (\nu_{w} (T^{n}))^{s}
\bigr ) \cr
&\qquad\qquad\qquad\qquad  \cdot | \det DT^n|_{E^s} |^{-1}\cdot
|\det DT^n |^{1-1/t} \, dLeb(w) \biggr )^{1/n}
\, . 
\end{split}
\end{equation}
\end{proposition}

\begin{proof}
The reader is invited to check  that there is $C_3>1$  (depending
on  $T$) so that
for  all $n\ge 1$, each cover $V$, all $i$, $j$, all $p\le 0$
$$ 
\max_{w\in \overline V_{ij,n}} (\lambda_{w} (T^n))^{p}
-\min_{w\in \overline V_{ij,n}}(\lambda_{w} (T^n))^{p}
\le n  C_3 (\lambda_{ij}(T^n))^p {\rm diam}\,  V \, ,
$$
and similarly for the $\nu_w$ (this is a bounded distortion argument). 
If  
$\ell_{ij}(n):=\max((\lambda_{ij}(T^n))^p, (\nu_{ij}(T^n))^s)=(\lambda_{ij}(T^n))^p
$ (the other case is similar) then
\begin{equation}\nonumber
\begin{split}
\max_{w\in \overline V_{ij}}(\ell_{ij}(n)& -\max( (\lambda_{w} (T^n))^{p}, (\nu_{w} (T^n))^{s}))
\le \ell_{ij} (n)-\min_{w\in \overline V_{ij}}(\lambda_{w} (T^n))^{p}\cr
&\quad\qquad\le  \max_{w\in \overline V_{ij}}(\lambda_{w} (T^n))^{p}
-\min_{w\in \overline V_{ij}}(\lambda_{w} (T^n))^{p} \, .
\end{split}
\end{equation}
Choose a partition $\XX=\cup_{i \in \II} W_i$
with $W_i\subset V_i$, and write $W_{ij}=W_i \cap T^{-n}W_j$.
Then 
\begin{equation}\nonumber
\begin{split}
&\sum_{i,j} \Leb(W_{ij}) \ell_{ij}(n) - \int_\XX 
 \max( (\lambda_{w} (T^n))^{p}, (\nu_{w} (T^n))^{s} ) \, d Leb\cr
&\qquad \le 
\sum_{i,j}   \Leb(W_{ij})\biggl ( \ell_{ij}(n) -  \min_{W_{ij}} 
 \max( (\lambda_{w} (T^n))^{p}, (\nu_{w} (T^n))^{s} )\biggr ) \cr
&\qquad \le
 C_3 n \cdot {\rm diam}\,  V \sum_{i, j} \Leb(W_{ij}) \ell_{ij}(n)  \, .  \cr
\end{split}
\end{equation}
Therefore, fixing $\delta \in (0,1)$, if  $V(n)$ satisfies   ${\rm diam}\,  V (n)= \delta/ (C_3n)$,  
$$
\sum_{i,j} \Leb(V_{ij,n}) \ell_{ij}(n) 
\le \frac{ \# V(n) }{(1-\delta)}  \int_\XX 
 \max( (\lambda_{w} (T^n))^{p}, (\nu_{w} (T^n))^{s} ) \, d Leb \, .
$$
Choose $V(n)$ and $\psi(n)$  with $\# V=O(n)$ and 
 $(\min_i \Leb (V_i) )^{-1}=O (n^{d})$,
ensuring that the derivatives of the $\psi_i$ 
from Lemma \ref{Leibniz} satisfy
$O(n^Q)$ bounds; for some $Q\ge 1$. 
Finally,  there is $C_4\ge 1$ so that 
$$
\frac{1}{\Leb(V_{ij,n})}  \le \frac{ C_4 }{\min_i \Leb (V_i))} \inf_{V_{ij,n}}
 |\det DT^n|_{E^u}|\, .
$$
for all $n$ and all covers $V$. Lemma ~\ref{comp} allows
to conclude, by a straightforward adaptation
of the proof of Theorem~\ref{mainthm}.
\end{proof}

\subsection{Proof of the Lasota-Yorke inequality}   \label{proofLY}

\begin{proof}[Proof of Lemma~\ref{comp}]
We replace $T$ by $T^n$  (the reader should keep in
mind that $A_{ij}$, $B_{ij}$, $D_{ij}$,
$\lambda_{ij}$, and $\nu_{ij}$ depend on $n$)
and drop the indices $i$, $j$. 
We study the action of the composition by $T$ on our symbol
$a_{p,q}(\xi,\eta)$
(see e.g. \cite[Chapter I.7, Proposition~7.1,
and Chapter I.8, Th\'eor\`eme 3]{AG}).

Taking a Taylor series 
of order $0$ (i.e., $k=1$ in the proof of  \cite[I.8, Lemme 4]{AG}), we find
that $(\Psi_{ij} a_{p,q}(\xi,\eta)^{Op} (f\circ T )) \circ T^{-1}$ decomposes as
\begin{equation}
\label{heart1}
\begin{split}
&( \Psi_{ij}\circ T^{-1})\cdot \bigl ( 
( a_{p,q}((DT)^{tr}_{T^{-1}(x,y)}(\xi,\eta))^{Op}( f)\cr
&\qquad \qquad\qquad+ r_1(x,y,\xi,\eta)^{Op} (f)
+r_2(x,y,\xi,\eta)^{Op} (f) \bigr )\, ,
\end{split}
\end{equation}
where $r_1$ and $r_2$ are described next.
The symbol $r_1(T(x,y),\xi,\eta)$ is a universal finite linear combination of
\begin{equation}\label{r1}
\begin{split}
&
 \int_{\real^d} du dv\, 
\int_{\real^{d}}  d\omega d\theta \, e^{-i(u,v)(\omega, \theta)}
\int_0^1 ds\, 
(1-s) s^j  \cr
&\cdot
\partial_{u_\ell} \bigl (e^{i(R_{(x,y)} (x+su,y+sv) )^{tr} (\xi,\eta)} 
\bigr )\cr
&\cdot
 (1+|s\omega+A_{(x,y)}\xi|^2+|s\theta+B_{(x,y)}\xi+D_{(x,y)}\eta|^2)^{p/2} \cr
&\cdot 
\partial_{\omega_\ell}
\bigl (1+|s\omega+ A_{(x,y)}\xi|^2 )^ {q /2}
 \chi  \biggl ( \frac{(s\omega+ A_{(x,y)} \xi,
s\theta+B_{(x,y)}\xi+D_{(x,y)}\eta)}
{1+|(A_{(0,0)}\xi,B_{(0,0)}\xi+D_{(0,0)}\eta)|}  \biggr ) \, ,
\end{split}
\end{equation}
where $j\in \{0, 1\}$, $1\le \ell \le d_s$, 
the function
$\chi:\real^{d}\to [0,1]$ is $C^\infty$ and compactly supported
in a suitable annulus, 
and 
$$R_{(x,y)}(u,v)=T(u,v)-T(x,y) -DT_{(x,y)} (u-x, v-y)\, .$$
To describe $r_2$, set $\tilde r_2=r_2\cdot (a_{p-1/2,q})^{-1}$,
so that $r_2^{Op}= \tilde r_2^{Op} a_{p-1/2,q}^{Op}$. 
(In the proof we shall use the notation
$\tilde r=r\cdot (a_{p-1/2,q})^{-1}$ several times.)
We claim that
$\tilde r_2^{Op}$ is a bounded operator on  $L^t( \real^d,\Leb)$,
so that 
\begin{equation}
\nonumber
\begin{split}
\|\Psi_{ij}\cdot \bigl ((r_2^{Op} f) \circ T\bigr) \|_{L^t} 
&\le C \sup_{V_{ijj}} |\det DT|^{-1/t}  \cdot
\biggl (\int | r_2^{Op}(f)|^t \, dLeb\biggr ) ^{1/t}\cr
&\le C_2(n) \|  f\|_{p-1/2, q, t,\real^d_j}\, .
\end{split}
\end{equation}
By \cite[Th\'eor\`eme~9]{CM} it 
suffices to show that for all $| \alpha|\le 1$ and all $ \beta$
we have $\sup | 
(1+|\xi|+|\eta|)^{|\beta| } \partial^{ \alpha}_{x,y} 
\partial^{\beta}_{\xi,\eta} \tilde r_2(x,y,\xi,\eta)|<\infty$.
This can be seen by observing that $r_2$ is made on the one
hand with contributions due to  $1-\chi$,
which have rapid decay in $1+|(\omega,\theta)|+
|(A_{(0,0)}\xi,B_{(0,0)}\xi+D_{(0,0)}\eta)|$
(by a small modification of \cite[p.58]{AG}, using
bounded distortion for $DT^n$). 
The other terms forming $r_2$ correspond
to a  $\partial_{\theta_\ell}$ derivative,
or to a  $\partial_{\omega_\ell}$,
but acting on a factor 
$
(1+|s\omega+A_{(x,y)}\xi|^2+|s\theta+B_{(x,y)}\xi+D_{(x,y)}\eta|^2)^{p/2}
$.
(Details are left to the reader, see \cite[p.60]{AG}.)

\smallskip
We may thus
concentrate on the first two terms in (\ref{heart1}).
The first one is called the {\it principal symbol.}

We get
$a_{p,q}((DT_{(x,y)})^{tr}(\xi,\eta))^{Op}=b^{Op}\circ  a^{Op}_{p,q}$
by setting $b(x,y,\xi,\eta)=a_{p,q}((DT_{(x,y)})^{tr}(\xi,\eta))/a_{p,q}(\xi,\eta)$.
Again by \cite[Th\'eor\`eme~9]{CM} it 
suffices to show that, up to replacing
$b$ by $b-r_3$, with $\tilde r_3^{Op}$ bounded
on each $L^t$, we have $\sup | (1+|\xi|+|\eta|)^{|\beta|} \partial^{ \alpha}_{x,y} 
\partial^{\beta}_{\xi,\eta} b|\le( C_\delta/2)
\max(\lambda^{p},\nu^{q+p})$, for all $| \alpha|\le 1$ and all $ \beta$.
Of course, we must also prove the same bounds for 
$r_1 \cdot (a_{p,q})^{-1}$ (modulo  
$r_4+r_5$, with $\tilde r_4^{Op}$ and $\tilde r_5^{Op}$
bounded on each $L^t(\real^d)$).

\medskip
Consider first $ \alpha= \beta=0$ and the principal symbol,
i.e., the bound for $\sup|b|$. 
For  $\xi\ne 0$ write $\nu_\xi= \sup_{x,y}|A_{(x,y)}\xi|/|\xi|$.  Then, setting
$\Lambda_1=\max_\xi ((1-2 \nu_\xi)/\nu_\xi)$,  we get for any  $|\xi|\ge \Lambda_1$
\begin{equation}\label{heart2}
(1+|A_{(x,y)}\xi|^2)^{q/2} 
\le (1+\nu_\xi^2|\xi|^2)^{q/2} \le  2 \nu_{\xi/|\xi|}^{q} (1+|\xi|^2)^{q/2}\, .
\end{equation}
For $|\xi|<\Lambda_1$  we always have
$
(1+|A_{(x,y)}\xi|^2)^{q /2} \le   (1+|\xi|^2)^{q/2} 
$.

\smallskip

If  $|\xi|\ge\max( \Lambda_1,|\eta|)$
then \footnote{Here we pay the price of $p< 0$.}
$$
(1+|A_{(x,y)}\xi|^2+|B_{(x,y)}\xi+D_{(x,y)}\eta|^2)^{p/2}
\le 2 \nu_{\xi/|\xi|}^{p} (1+|\xi|^2+ |\eta|^2)^{p/2}\, .
$$

For  $\eta\ne 0$ write $\lambda_\eta= \inf_{x,y} |D_{(x,y)}\eta|/|\eta|$. 
Fix $\Lambda_2=\max_\eta \lambda_\eta$.
If $\epsilon <1/4$ and $|\eta|\ge \max(\Lambda_2, |\xi|)$ then
\footnote{If $|\xi|$ is small but $|\eta|$ is large we need $p< 0$ to get a contraction here.}
$$
(1+|A_{(x,y)}\xi|^2+|B_{(x,y)}\xi+D_{(x,y)}\eta|^2)^{p/2}
\le 3 \lambda^{p} (1+|\xi|^2+|\eta|^2)^{p/2}\, .
$$

Finally, if $|\xi|\le \Lambda_1$, and $|\eta|\le \Lambda_2$, there is
$C_{\Lambda_1, \Lambda_2}(n)$  with
\begin{equation}\label{heartn}
(1+|A_{(x,y)}\xi|^2+|B_{(x,y)}\xi+D_{(x,y)}\eta|^2)^{p/2} \le C_{\Lambda_1, \Lambda_2}
 (1+|\xi|^2+|\eta|^2)^{p/2-1}\, .
\end{equation}
We  include
the above contribution in  $r_3$.

\medskip
To bound $\sup_{x,y,\xi,\eta}|r_1\cdot (a_{p,q})^{-1}|$,  
multiply the integrand of (\ref{r1})
by 
$$
\biggl (1- \tilde \chi \bigl (\frac{ s\omega+ A_{(x,y)} \xi }{1+|A_{(0,0)}\xi|}
\bigr )  \biggr ) + 
\tilde \chi  \bigl  ( \frac{s\omega+ A_{(x,y)} \xi}{1+|A_{(0,0)}\xi|} \bigr )\, , $$ 
where $\tilde \chi:
\real^{d_s} \to [0,1]$ is $C^\infty$ and compactly supported
in an annulus.
We consider separately the two terms in this decomposition:

\smallskip
The term containing  $\chi \cdot (1-\tilde \chi)$ enjoys
$C_k(n)(1+|A_{(0,0)}\xi|+|\omega|)^{-k}$ rapid decay (adapting \cite[p.58]{AG}).
By choosing  first $k$ and then $\Lambda_3$ we get a bound $(C_\delta/4)
\max(\lambda^{p},\nu^{q+p})$  for $|\xi|\ge \Lambda_3$.
If $|\xi|\le \Lambda_3$, we use
that if $|\eta|>\max( |\xi|,\Lambda_4)$ then
\begin{equation}\label{suivantes}
\begin{split}
\sup_{s,(\omega,\theta),(u,v),(x,y)}&
(1+|s\omega+A_{(x,y)}\xi|^2+|s\theta+B_{(x,y)}\xi+D_{(x,y)}\eta|^2)^{p/2}\cr
 &\quad \cdot \chi  \biggl ( \frac{(s\omega+ A_{(x,y)} \xi,
s\theta+B_{(x,y)}\xi+D_{(x,y)}\eta)}
{1+|(A_{(0,0)}\xi,B_{(0,0)}\xi+D_{(0,0)}\eta)|} \biggr )\cr
&\qquad \le 2 \lambda^p (1+|\xi|^2+|\eta|^2)^{p/2} \, .
\end{split}
\end{equation}
The compact set $\{ |\xi|\le \Lambda_3, |\eta|\le \Lambda_4\}$ gives
rise to a term $r_4$.

For  the $\chi \cdot \tilde \chi$ term,
 use the  ideas exploited for the
principal symbol (see also -- again -- \cite[p.60]{AG}) to get  a bound  $(C_\delta/4)
\max(\lambda^{p},\nu^{q+p})$, up to a  perturbation $r_5$.
In particular, if $|\xi|$ is large with respect to $|\eta|$  then
\begin{equation}\label{term}
\begin{split}
\sup_{s,(\omega,\theta),(x,y)} &(1+|s\omega+A_{(x,y)}\xi|^2+|s\theta+B_{(x,y)}\xi+D_{(x,y)}\eta|^2)^{p/2}\cr  
&\qquad \cdot \partial_{\omega_\ell}
\bigl ( (1+|s\omega+ A_{(x,y)}\xi|^2 )^ {q /2}
 \chi  (\cdots) \cdot \tilde \chi (\cdots) \bigr ) 
\end{split}
\end{equation}
is bounded by
$C (1+|\xi| +|\eta|)^{p-1} (1+|\xi| )^ q$
(giving a  contribution $r_5$);
while if $|\eta|$ is large then 
(\ref{term}) is bounded by
$2 \lambda^{p} (1+|\xi|^2 +|\eta|^2)^{p/2} (1+|\xi|^2 )^{q/2}$.

\smallskip

The control of the derivatives of
$b$ and $r_1\cdot( a_{p,q})^{-1}$, i.e., the case of nonzero $| \alpha|+| \beta|$,
is straightforward although rather tedious.
\end{proof}

\subsection{The essential spectral radius of $\MM$}
\label{last}

Let $T$ be a $C^\infty$ Anosov diffeomorphism  on a compact
manifold, with a $C^\infty$ unstable foliation. 
Let $W^{p,q,t}(\XX,T^{-1})$ denote the Banach space in
Definition \ref{defnorm}.
(We use now stable foliation of $T^{-1}$, i.e. the unstable foliation of $T$,
in other words,
$W^{p,q,t}(\XX,T^{-1})=(1+\Delta_u)^{-q/2} (1+\Delta)^{-p/2} ( L^t(\XX))$.)

\begin{theorem}[Essential spectral radius of $\MM$]\label{mainthm2}
For any $p<0$ and $s>0$, the essential spectral radius of $\MM$  on
$W^{p,s-p,t}(\XX,T^{-1})$  is not larger than
$\lim_{n \to \infty} \sup_\XX |\det DT^n|^{-(t-1)/tn}  \cdot \rho^{(-s,-p)}_\infty(T)$
for   $t \in(1,\infty)$.
\end{theorem}

\begin{proposition}\label{L1,2}
For any $p<0$, $s>0$, and   $1< t< \infty$,
the essential spectral radius of $\MM$  on
$W^{p,s-p,t}(\XX,T^{ -1})$  is  
\begin{equation}
\begin{split}
&\le 
\lim_{n\to \infty}
\biggl ( \int_\XX
\max\bigl ( (\lambda_{w} (T^n))^{-s}, (\nu_{w} (T^n))^{-p}
\bigr )\cr
& \qquad\qquad\qquad\qquad \cdot | \det DT^{n}|_{E^u} |\cdot
|\det DT^{n} |^{-(t-1)/t} \, dLeb(w) \biggr )^{1/n} 
\, . 
\end{split}
\end{equation}
\end{proposition}

\begin{proof}[Proof of the theorem and the proposition]
Adapt the proofs of Theorem \ref{mainthm} and
Proposition~\ref{L1}, using distortion estimates to
bound (\ref{opM}) when exploiting Lemma~\ref{Leibniz}
for a weight $(1/|\det DT^n |)\circ \kappa_i^{-1}$
(see also the comments before  Corollary~\ref{compint}).
\end{proof}

\appendix

\section{Operators $\MM_t$ and $\LL_t$}

\begin{theorem}\label{mainthm3}

\noindent\phantom{nimportequoi}

\noindent (1) If the stable foliation of $T$ is
$C^\infty$ then for any $p<0$ and $s>0$ with $q=s-p$ integer,
there exists $t_1(q)>1$ so that
the essential spectral radius of $\LL_t$ on
$W^{p,s-p,t}(\XX)$  is 
 $\le \rho^{(p,s)}_\infty(T)$, for each  $1<t<t_1$.

\noindent (2) If the unstable foliation of $T$ is
$C^\infty$ then for any $p<0$ and $s>0$ with $q=s-p$ integer 
there exists $t_2(q)\ge 1$ so that
the essential spectral radius of $\MM_t$  on
$W^{p,s-p,t}(\XX,T^{-1})$  is 
 $\le \rho^{(-s,-p)}_\infty(T)$, for each  $t_2 < t <\infty$.
\end{theorem}

For any integer
$q\ge 1$, there exist $C\ge 1$ and $t_1(q) >1$ so that for 
every  $\gamma'$ with $1\le |\gamma'|\le q$,
all $1\le t < t_1$
and all $n\ge 1$, setting $h(x,y)=|\det DT^{n}| ^{1/t} \circ \kappa_i^{-1}(x,y)$,
then
$|\partial^{\gamma'}_x h(x,y)|\le C h(x,y)$.
(If $q=1$ we may take $t_1=\infty$.)
Theorem \ref{mainthm3}  is therefore a consequence of the following corollary of the 
proof of Lemma~\ref{comp} combined
with a  refinement of the Leibniz
formula for $a_{p,q}$ if $q\in \integer_+$, Lemma~\ref{Leibnizint}. 

\begin{corollary}[More on Lasota-Yorke]\label{compint}
There exist $\delta_0 >0$ and $C_0$ so that, for all $V$ with
${\rm diam} \, V <\delta_0$
and   $n\ge 1$,  there exists $C(n)>1$ so that
for any multi-index $\gamma'$ with $|\gamma'|\le q$, any $f \in W^{p,q,t}$,
compactly supported in $U_j$,
and each $C^\infty$ function
$\Psi_{ij}:\real^d_i\to [0,1]$ compactly supported in $U_{ij,n}$
\begin{equation}\nonumber
\begin{split}
\|   \Psi_{ij}\cdot & |\det DT^n \circ \kappa^{-1}_i|^{1/t} \cdot a_{p,0}\, \partial_x^{\gamma'} (f \circ   T_{ij}^n)\|_{L^t} \cr
&\qquad\le C_0 
\max((\lambda_{ij}(T^n))^{p},(\nu_{ij}(T^n))^{q+p})  
\|  f \|_{p,q,t,\real^d_j}\cr
&\qquad\qquad\quad+ C(n) \|  f \|_{p-1/2,q,t,\real^d_j} \, , \, \forall
p\le 0\, , \, q\ge -p\, , \, 1<t<\infty \, .
\end{split}
\end{equation}
\end{corollary}

\begin{lemma}[Leibniz formula for integer derivatives]\label{Leibnizint}
Let $1<t<\infty$. Let $q \in \integer_+$ and let $p \in \real$.
There exists $C\ge 1$, and for every compactly supported
 $h\in C^\infty(I^d)$ there exists $C(h)>0$ so that for each
$f \in W^{p,q,t}(\real^d)$, we have
$a^{Op}_{p,q}( h \cdot f) =h \cdot a^{Op}_{p,q}(f)+ g_1 + g_2$
with 
$$
\|g_1\|_{L^t} \le C 
\sum_{\stackrel{|\gamma'| =q-1}{ \gamma'_1+\gamma_2'=\gamma'}}
 \| \partial^{\gamma'_1}_x h \cdot a_{p, 0}^{Op}( \partial^{\gamma'_2}_x f)\|_{L^t}\, , \,\,
\|g_2\|_{L^t} \le C(h)  \| f\|_{p-1,q,t} \, . 
$$
\end{lemma}

\begin{proof} 
Decompose $a_{p,q}^{Op}=a_{p,0}^{Op} \circ a_{0,q}^{Op}$.
The proof of Lemma~\ref{Leibniz} gives that if $\tilde h \in
C^\infty(I^d)$ is compactly supported, then there is $C(\tilde h)\ge $ so that
for all 
$\tilde f \in W^{p,0,t}(\real^d)$,   we have
$
a^{Op}_{p,0}( \tilde h \cdot \tilde f) =\tilde h \cdot a^{Op}_{p,0}(\tilde  f)+
\tilde g
$
with $\| \tilde g \|_{L^t}\le C(\tilde h) \|\tilde f \|_{p-1, 0, t}$.

Since $q$ is an integer,
$a^{Op}_{0,q}(h \cdot f)$ decomposes as:
\begin{equation}
\nonumber
a^{Op}_{0,q}(h \cdot f)=
\begin{cases} h\cdot f+ 
\sum_{j=1}^\ell
\binom {\ell}{j} \bigl (\partial^2_{x_1}+ \cdots +\partial^2_{x_{d_s}} \bigr )^j (h \cdot f)&\cr
\qquad \qquad\qquad\qquad \qquad\qquad\qquad\text{if } q=2\ell \text{ is even,}\cr
\mu_1 * a^{Op}_{0,q-1}(h \cdot f) +
\mu_2 * \bigl ( \sum_{j=1}^{d_s} R_{x_j} ( \partial_{x_j} a_{0, q-1}^{Op}(h\cdot f )\bigr )\cr
\qquad\qquad\qquad\qquad \qquad\qquad\qquad \text{if } q=2\ell+1 \text{ is odd,}\cr
\end{cases}
\end{equation}
where $\mu_1$ and $\mu_2$ are finite measures (which do not depend on
$h$ or $f$) and $*$ denotes convolution. Indeed,
$q=2\ell$ is even, just recall that
$a_{0,q}^{Op}=(1+\Delta_s)^{q/2}= (1+\sum_{j=1}^{d_s} \partial^2_{x_j})^{\ell}$.
If $q=2\ell+1$, recall \cite[V.3.2--V.3.4]{St} that
$
(1+\Delta_s)^{1/2} (\varphi)
=  \mu_{1} * \varphi +  \mu_2 * 
\bigl (  \sum_{j=1}^{d_s} R_{x_j} ( \partial_{x_j} \varphi ) \bigr )
$,
where $R_{x_j}$ is the Riesz transform
\cite[III.1]{St}.
Finally, use the ordinary Leibniz formula
for partial derivatives, that $a^{Op}_{p,0}$
commutes with each $R_{x_j}$, that $R_{x_j}$ is bounded on  $L^t$,
and  that if $\mu$ is a measure, with total mass $|\mu|$, then 
$a_{p,0} (\mu * f)= \mu * a_{p,0} (f)$ and
$\|\mu * f \|_{L^t}\le |\mu| \cdot \|f\|_{L^t}$.
\end{proof}

\bibliographystyle{amsplain}

\end{document}